\documentclass[12pt]{amsart}
\usepackage{graphicx,amssymb}

\theoremstyle{definition}

\theoremstyle{remark}

\numberwithin{equation}{section}

\newcommand{\R}{\mathbb R}
\newcommand{\C}{\mathbb C}

\def\TagOnRight

\def\R {\mathbb{R}}

\newcommand{\be}{\begin{equation}}
\newcommand{\ee}{\end{equation}}
\newcommand{\bea}{\begin{eqnarray}}
\newcommand{\eea}{\end{eqnarray}}
\newcommand{\Bea}{\begin{eqnarray*}}
\newcommand{\Eea}{\end{eqnarray*}}
\newcommand{\bt}{\begin{Theorem}}
\newcommand{\et}{\end{Theorem}}

\newcommand{\bpr}{\begin{Proposition}}
\newcommand{\epr}{\end{Proposition}}

\newcommand{\bl}{\begin{Lemma}}
\newcommand{\el}{\end{Lemma}}
\newcommand{\bi}{\begin{itemize}}
\newcommand{\ei}{\end{itemize}}

\newtheorem{Definition}{Definition}[section]
\newtheorem{Theorem}[Definition]{Theorem}
\newtheorem{Lemma}[Definition]{Lemma}
\newtheorem{Proposition}[Definition]{Proposition}

\begin{document}
\baselineskip16pt

\title[]{ Benedick's theorem for the Heisenberg Group}
\author{E.K. Narayanan and P. K. Ratnakumar}%
\address{Department of Mathematics, Indian Institute of Science,
Bangalore -12, India}
\address {Harish-Chandra Research Institute, Allahabad-211019,
 India.}
 \email {naru@math.iisc.ernet.in, ratnapk@hri.res.in}
\thanks{The first author was supported in part by a grant from
UGC via DSA-SAP.}%

\begin{abstract}
If an integrable function $f$ on the Heisenberg group is
supported on the set $B \times \R$ where
$B \subset \C^n$ is compact and the group Fourier transform
$\hat{f}(\lambda)$ is a finite rank operator for all $\lambda
\in \R \setminus \{0\},$ then $f \equiv 0.$
\end{abstract}

\maketitle

\section{Introduction}
\noindent

The uncertainty principle says that a nonzero function and its
Fourier transform cannot both be sharply localized. There are
several manifestations of this principle. We refer the reader
to the excellent survey article by Folland and Sitaram [6]
and also the monograph by S. Thangavelu [9].

In this paper we are interested in a variant of
Benedick's theorem on
the Heisenberg group. Recall that the Benedick's theorem [2]
states the following. Let $f \in L^2(\R^n),$ if both the sets
$\{x \in \R^n : f(x) \neq 0 \}$ and
$\{\xi \in \R^n : \hat{f}(\xi) \neq 0 \}$
have finite Lebesgue measure, then $f \equiv 0$. In the context
of non commutative Lie groups the Fourier transform is an
operator valued function. We measure
the ``smallness'' of the Fourier transform in terms of the
rank of these operators.

To state our result, we need to recall briefly the
representation theory of the Heisenberg group.
The Heisenberg group $\mathbb{H}^n$ is topologically
$\C^n \times \R$, with the group law
$$(z,t)\cdot (w,s)= (z+w, t+s+ \frac{1}{2}\Im (z \cdot \bar{w})).$$
 Under this group law, $\mathbb{H}^n$ becomes a two step nilpotent
Lie group with center $Z= \{ 0\} \times \R.$ The infinite
dimensional irreducible unitary representations
of $\mathbb{H}^n$ are parametrized by
 $\lambda \in \R \setminus \{0\}$. Each such $\lambda$ defines a
representation $\pi_\lambda,$ realized on $ L^2(\R^n)$
by $$\pi_\lambda(z,t)\varphi(\xi) =e^{i \lambda t} \,
 e^{i \lambda(x \cdot \xi + \frac{1}{2}x \cdot y)}
\varphi(\xi + y)$$ where $z= x+ i y$ and $\varphi \in L^2(\R^n).$
The representation $\pi_\lambda$ is clearly unitary and
it is well known that they are irreducible on $L^2(\R^n)$.
In fact, a famous theorem of Stone and von Neumann says that
any irreducible unitary representation of $\mathbb{H}^n$, that is non
 trivial at the center is (unitarily) equivalent to $\pi_\lambda$
for some $\lambda$ (see [5]).

If $f \in L^1(\mathbb{H}^n)$, we can define the group Fourier transform
 by $$\hat{f}(\lambda)= \int_{\mathbb{H}^n}f(z,t) \pi_\lambda(z,t)\, dz \, dt.$$
Since $\pi_\lambda$ is an isometry, a simple norm estimate
shows that $\hat{f}(\lambda)$ is a bounded operator on $L^2(\R^n)$.
Moreover, if $f\in L^2(\mathbb{H}^n)$ then $\hat{f}(\lambda)$ turns out to be a Hilbert
Schmidt operator and the Plancherel theorem for the Heisenberg group reads as
$$\int_{\mathbb{H}^n}|f(z,t)|^2 dz \, dt = (2 \pi)^{-n-1}
 \int \| \hat{f}(\lambda)\|_{HS}^2 |\lambda|^n d\lambda.$$
Our main result is the following :

\bt Let $f \in L^1(\mathbb{H}^n)$ is supported on  a set of the form
$B \times \R$, where $B \subset \C^n.$

I. If $B$ is a compact set and   $\hat{f}(\lambda) $
 is a finite rank operator for all $\lambda \neq 0$,
then $f\equiv 0$.

II. If $B$ has finite Lebesgue measure and $\hat{f}(\lambda) $
 is a rank one operator for all $\lambda \neq 0$, then $f\equiv 0$.
\et

\noindent
{\bf Remark I:} Note that our result is in
sharp contrast with the situation on other Lie groups.
For example, in the Euclidean case, the Fourier transform of any nontrivial $f \in L^1(\mathbb{R}^n)$
gives rise to a rank one operator on $L^2(\mathbb{R}^n)$; via multiplication by $\hat{f}(\xi)$.
 Next, if $G$ is a non compact connected semisimple Lie group
 and $K$ is a maximal compact subgroup of $G,$
then it can be shown that a function in $L^1(G/K)$
 has a Fourier transform, which is a rank one operator.
 More generally, if $f \in L^1(G)$ transforms
according to a fixed unitary irreducible representation of
the compact group $K$ on the right, then the group
Fourier transform of $f$ is a finite rank operator.

\noindent
{\bf Remark II:} In [1] the authors study ``Qualitative 
uncertainty Principle'' for unimodular groups. Let $G$ be such
a group and $\hat{G}$ be its unitary dual. Let $m$ denote the
Haar measure on $G$ and $\hat{m},$ the Plancherel measure on
$\hat{G}.$ One of the results in [1] states that, if $\{
x \in G:~f(x) \neq 0 \} < m(G)$ and $\int_{\hat{G}}
rank (\pi(f))~d\hat{m} < \infty$ then $f \equiv 0.$  When
$G$ is the Heisenberg group, the above conditions will force
the Fourier transform to be supported on a set of finite 
(Plancherel) measure in addition to the finite rank condition.
Notice that, Theorem 1.1 requires only the finite rank condition.
We thank Michael Lacey for pointing out this reference.
We also refer the reader to [7] for a Benedick's type theorem
on the Heisenberg group which is mainly a $t$-variable theorem.

In the rest of this section, we recall the necessary
details about the Weyl transform and the
Fourier-Wigner transform.
For a suitable function $g$ defined on $\C^n$, the
 $\lambda-$Weyl transform is defined to be the
 operator $$W_\lambda(g) = \int_{\C^n} g(z)\,
\pi_\lambda(z) \, dz$$ where $\pi_\lambda(z)
= \pi_\lambda(z,0).$ Clearly  $W_\lambda(g)$ defines a bounded
operator on $L^2(\R^n)$, if $g \in L^1(\C^n)$. 
For $g \in L^2(\C^n),~ W_\lambda(g)$ is a Hilbert-Schmidt operator
 and we have the Plancherel Theorem [8]:
$$\int_{\C^n}|g(z)|^2 dz= (2 \pi \, |\lambda|)^{-n} \, \| W_{\lambda}
 (g)\|_{HS}^2 .$$

The $\lambda-$twisted convolution of two
functions $F$ and $G$ on $\C^n$ is defined
 to be $$F\times_\lambda  G(z) = \int_{\C^n}
 F(z-w) G(w) e^{\frac{i \lambda}{2} \Im (z \cdot \bar{w})} dw.$$
It is known that $W_\lambda (F\times _\lambda G)
 = W_\lambda (F) W_\lambda ( G)$.
When $\lambda =1$, we write $F \times G $
instead of $F \times_1 G$ and call
it the twisted convolution of $F$ and $G.$ Similarly
$W_1(F)$ will be denoted by $W(F)$, and called
the Weyl transform of $F.$

Let $\phi_1$ and $\phi_2$ belong to $L^2(\R^n).$ The
Fourier-Wigner transform
of $\phi_1$ and $\phi_2$ is a function on $\C^n$ and is defined
by $$A(\phi_1, \phi_2)(z) = \langle \pi_1(z)\phi_1, \phi_2
\rangle .$$ The Fourier-Wigner transform satisfies the `orthogonality relation',
\bea
\int_{\C^n}~A(\phi_1,\phi_2)(z)~\overline{A(\psi_1, \psi_2)
(z)}~dz = (2 \pi)^n \, \langle \phi_1, \psi_1 \rangle~\langle \psi_2, \phi_2
\rangle .
\eea
In fact, if $\{\phi_i : i \in \mathbb{N} \}$ is an orthonormal basis for $L^2(\R^n),$ then the collection
$\{ A(\phi_i, \phi_j) : i, j \in \mathbb{N} \} $ form an orthonormal basis for
$L^2(\C^n),$ see [8]. In particular, if $F \in L^2(\C^n)$ is orthogonal to
$A(\phi, \psi)$ for all $\phi, \psi \in L^2(\R^n)$ then $F\equiv 0$.

We finish this section with the following theorem (see [3] or [4])
which will be used later.

\bt
Let $F(z) = A(\phi_1, \phi_2)(z)$ where $\phi_1, \phi_2 \in
L^2(\R^n).$ If the set $\{ z:~F(z) \neq 0 \}$ has finite
Lebesgue measure then $F \equiv 0.$
\et

\section{Proof of The main result}

We start with the following lemma:

\bl
Let $h_j \in L^2(\R^n), ~ j=1,2,...,N$,
and set, for $y \in \R^n$  $$K_y(\zeta) =
\sum_{j=1}^N {\overline{h_j(\zeta)}} h_j(\zeta+y).$$
If $K_y(\zeta) =0$ for almost all $\zeta \in \R^n$
 and $|y| \geq R$, then each $h_j$ is compactly supported.
\el

\noindent
{\bf {Proof:} } Since each $h_j \in L^2(\R^n)$,
 there exist a set $A$ of Lebesgue measure zero
such that $|h_j(\zeta)|< \infty $ for every
$\zeta \in \R^n \setminus A,$ for $j=1,2\dots,N.$

We work with a fixed representative $h_j$ for each of the class $[h_j] \in L^2(\R^n)$
for which pointwise evaluation makes sense. Hence, for $\zeta \in \R^n \setminus A$,
 $$H(\zeta) = (h_1 (\zeta), h_2(\zeta), \dots h_N(\zeta)) \in \mathbb{C}^N.$$
If $h_j$ are non zero, choose $\zeta_1 \in \R^n \setminus A$
so that $H(\zeta_1)$ is a non zero vector. Let $B_R(\zeta_1)$
be the open ball of radius $R$ centered at $\zeta_1.$ If there
is no $\zeta \in \R^n \setminus (B_R(\zeta_1) \cup A)$
such that $H(\zeta)$ is a non zero vector, we are done.
Otherwise, choose $\zeta_2 \in \R^n \setminus (B_R(\zeta_1) \cup
A)$ so that $H(\zeta_2)$ is non zero. By the hypothesis,
$H(\zeta_1)$ and $H(\zeta_2)$ are orthogonal vectors in $\C^N.$
We repeat this process. That is, if $j \leq N,$ choose $\zeta_j
\in \R^n \setminus (\cup_{l=1}^{j-1} B_{R}(\zeta_l) \cup A)$
such that $H(\zeta_j)$ is a non zero vector in $\C^N$ (if there
is no such $\zeta_j$ we are done). By the hypothesis $H(\zeta_j)$
are orthogonal to each other for $j = 1, 2, \cdots N .$
Now, if $\zeta \in \R^n
\setminus (\cup_{j=1}^N B_R(\zeta_j) \cup A)$ then $H(\zeta)$ is
orthogonal to $H(\zeta_j)$ for all $j = 1, 2, \cdots N .$ It
follows that $H(\zeta)$ is zero for $\zeta \in \R^n \setminus
(\cup_{j=1}^N B_R(\zeta_j) \cup A)$ which finishes the
proof.~\hfill $\Box$

Our next result is a Benedick's type theorem for the Weyl transform
and is a crucial step in the proof of the main theorem.

\bt
Let $F \in L^1(\C^n)$ be compactly
supported. If the Weyl transform $W(F)$ of
$F$ is a finite rank operator, then $F \equiv 0.$
\et

{{\bf Proof:}}
 Let ${\overline{G(z)}} = F^* \times F (z)$, where $F^*(z) =
\overline{F(-z)}$. Then $\bar{G}$ is compactly supported and
$\bar{G}=0$ if and only if $F\equiv 0$, by the Plancherel theorem for the Weyl transform.
Now $W(\bar{G})=W(F)^* W(F)$ is a finite rank, positive,
Hilbert-Schmidt operator and hence by the spectral theorem, we have
\bea W(\bar{G}) \phi = \sum_{j=1}^N \lambda_j \langle \phi,
\phi_j \rangle \phi_j \eea
where $\{\phi_1,..., \phi_N \}$ is an orthonormal basis for the range of $W(\bar{G}),$
with $W(\bar{G})\phi_j=\lambda_j\phi_j$ and $\lambda_j\geq 0$. Hence
 \bea
\langle W(\bar{G})\phi, \psi \rangle =
 \sum_{j=1}^N \lambda_j \langle \phi, \phi_j
\rangle \langle \phi_j , \psi \rangle.
\eea By (1.1), the above equals
\bea
(2 \pi)^{-n}\, \sum_{j=1}^N~\lambda_j~\int_{\C^n}~A(\phi, \psi)(z)~\overline{
A(\phi_j, \phi_j)(z)}~dz
\eea
Also by the definition of the Weyl transform,
\bea \langle W(\bar{G})\phi,
\psi \rangle = \int_{\C^n} \bar{G} (z)
\, A( \phi, \psi)(z) \, dz. \eea
From (2.3) and (2.4) it follows that
$$G(z)= \sum_{j=1}^N
A( h_j, h_j)(z),$$
where $h_j(z) = (2 \pi)^{-\frac{n}{2}}\sqrt{\lambda_j} \, \phi_j(z).$
 Writing $G_y(x)= G(z) $ for $z= (x+iy)$, the above identity reads as
\bea G_y(x) =
\int_{\R^n}e^{i(x \cdot \zeta + \frac{1}{2} x \cdot y )}
\left( \sum_{j=1}^N h_j (\zeta+y) {\overline{h_j(\zeta)}}
\right) d \zeta .
\eea

Since $G$ is compactly supported,
there exists $R>0$ such that, $G_y \equiv 0$ if $|y|\geq R$. Then (2.5)
implies that $\sum_{j=1}^N h_j (\zeta+y) \overline{h_j(\zeta)}=0$
for almost every $\zeta \in \R^n$, provided $|y| \geq R$.

Lemma 2.1 now implies that each $h_j$ is compactly supported
and hence $ \sum_{j=1}^N h_j (\zeta+y) \overline{h_j(\zeta)}$ is also compactly supported
in $\zeta$ for each $y\in \mathbb{R}^n$. In view of (2.5), we conclude that $G_y \equiv 0$ 
for each $y \in \R^n$, hence the proof. \hfill $\Box$\\

Now, we are in a position to complete the proof of the main theorem.\\

{{\bf Proof of Theorem 1.1 :} }
Let $f^\lambda (z)$  denote the partial Fourier transform of $f$
in the $t-$variable. That is $$f^\lambda (z) = \int_{\R}~f(z, t)~
e^{i\lambda t}~dt.$$ Then a simple computation shows that
$\hat{f}(\lambda) =W_\lambda(f^\lambda).$

We start with the proof of $(II)$ in Theorem 1.1. By the hypothesis
we have that $f^\lambda(z)$ is supported in the set $B$ (which
has finite Lebesgue measure) and $\hat{f}(\lambda) = W_\lambda(
f^\lambda)$ is a rank one operator for all $\lambda.$
 We will assume that $\lambda=1$ and prove that
$f^\lambda\equiv 0$. The general case is no different.

It suffices to
show that, if $F\in L^1(\C^n)$ is supported on a set of finite
measure and $W(F)$ is a rank one operator then $F\equiv 0$. This
immediately follows from Theorem 1.2 once we show that $F$ is
the Fourier-Wigner transform of two functions in $L^2(\R^n)$.
 For this, let $\bar{G}=F$. Since $W(\bar{G})$ is a rank one
operator, we
have $\psi_1,\psi_2 \in L^2(\R^n)$ such that $$W(\bar{G})\varphi=
\langle \varphi, \psi_1 \rangle~
 \psi_2,  ~\forall \varphi \in L^2(\R).$$
Hence, if $\psi \in L^2(\R^n)$ we have
\Bea
\langle W(\bar{G})\varphi, \psi \rangle &=&
\int_{\C^n} \bar{G}(z) \, \langle \pi_1(z)\varphi, \psi
\rangle \, dz\\
&=& \langle \varphi, \psi_1 \rangle \, \langle  \psi_2, \psi
\rangle \\
&=& (2\pi)^{-n}~\int_{\C^n}~A(\varphi, \psi)(z)~\overline{
A(\psi_1, \psi_2)(z)}~dz, \\
\Eea
where the last step follows from (1.1). It follows that
$G(z)= A(\psi_1, \psi_2)(z).$ 

To prove (I), we proceed as above. Taking the Fourier
transform in the $t-$ variable reduces the problem to
$\C^n$. As above we assume that $\lambda=1.$ It suffices to
show that if $F \in L^1(\C^n)$ is compactly supported and
$W(F)$ is of finite rank, then $F\equiv 0$. But this
is precisely the content of Theorem 2.2, hence finishes the proof. \hfill $\Box$


\begin{thebibliography}{99}

\bibitem {AL} D. Arnal and J. Ludwig, QUP and Paley-Wiener
properties of unimodular, especially nilpotent, Lie groups,
{\em Proc. Amer. Math. Soc.} 125 (1997), no.4, 1071-1080. \\

\bibitem {B} M. Benedicks, On Fourier transforms
of functions supported on a set of finite Lebesgue measure,
{\em J. Math. Anal.Appl.}, 106 (1985), 180-183 .\\

\bibitem {J} P. Jamming, Principe d'incertitude qualitatif
et reconstruction de phase pour la transform de Wigner,
{\em C. R. Acad. Sci. Paris, I Math}. 327 (1998), no.3, 249-254. \\

\bibitem {Ja} A. J. E. M. Janssen, Proof of a conjecture
on the supports of Wigner distributions, {\em J. Fourier Anal.
Appl}. 4 (1998), no.6, 723-726. \\

\bibitem {} G. B. Folland, {\em Harmonic Analysis in Phase
space}, Annals of math studies 122, Princeton Univ. Press.,
Princeton NJ, (1989).\\

\bibitem {FS} G. B. Folland and A. Sitaram,
Uncertainty principle, a mathematical survey,
{\em Journal of Fourier Anal. Appl}. 3(1997), no.3, 207-238.\\

\bibitem {PS} J. F. Price and A. Sitaram, Functions and their
Fourier transforms with supports of finite measure for
certain locally compact groups, {\em J. Funct. Anal.} 79 (1988),
no.1, 166-182.

\bibitem {T1} S. Thangavelu , {\em Lectures on Hermite and Laguerre
expansions,} Mathematical notes, 42, Princeton Univ.
press, Princeton (1993).\\

\bibitem {T2} S. Thangavelu , {\em An introduction to the
uncertainty principle}, Progress in Math. Vol. 217,
Birkhauser Boston, Inc.,Boston, MA(2004).\\


\end{thebibliography}
\end{document}